\documentclass[a4paper,11pt]{amsart}
\usepackage{amssymb, graphics,setspace}
\input epsf
\newtheorem {thm}   {Theorem}
\newtheorem {lemma} {Lemma} [section]
\newtheorem {prop}[lemma]{Proposition}

\theoremstyle{definition}

\newtheorem{rem}[lemma]{Remark}

\newtheorem {cor}[lemma]{Corollary}
\newtheorem{defn}[lemma]{Definition}

\newcounter{numl}

\newcommand{\labelnuml}{\textup{(\roman{numl})}}
\newenvironment{numlist}{\begin{list}{\labelnuml}%
{\usecounter{numl}\setlength{\leftmargin}{0pt}%
\setlength{\itemindent}{2\parindent}%
\setlength{\itemsep}{\smallskipamount}\def
\makelabel ##1{\hss \llap {\upshape ##1}}}}{\end{list}}

\newcommand{\R}{{\mathbb R}}
\newcommand{\C}{{\mathbb C}}

\newcommand{\cF}{{\mathcal F}}
\newcommand{\cA}{{\mathcal A}}

\newcommand{\cO}{{\mathcal O}}

\newcommand{\Scal}{\mathit{Scal}}

\newcommand{\grad}{\mathop{\mathrm{grad}}\nolimits}

\newcommand{\symprod}{\mathbin{\raise1pt\hbox{$\scriptstyle\bigcirc$}}}
\newcommand{\ra}{\rightarrow}

\newcommand{\Mpc}{p_{\mathrm{c}}}

\newcommand{\Ref}[1]{(\ref{#1})}

\makeatletter
\renewcommand\@biblabel[1]{#1.}     % Arabic numbers, no brackets
\makeatother
\begin{document}

%\doublespace
\title{Generalized Quasi-Einstein Metrics on Admissible Manifolds}

%\date{\today}
\begin{abstract} We prove that an admissible manifold (as defined by
Apostolov,  Calderbank, Gauduchon and T{\o}nnesen-Friedman), arising from a
base with a local K\"ahler product of constant scalar curvature metrics,
admits Generalized Quasi-Einstein K\"ahler metrics (as defined by D. Guan) in
all ``sufficiently small'' admissible K\"ahler classes. We give an
example where the existence of Generalized Quasi-Einstein metrics fails
in some K\"ahler classes while not in others. We also prove
an analogous existence theorem for an additional metric type, defined by
the requirement that the scalar curvature is an affine combination of
a Killing potential and its Laplacian.
\end{abstract}
\author[G. Maschler]{Gideon Maschler}
\address{Department of Mathematics and Computer Science
\\Clark University\\Worcester\\Massachusetts 01610\\USA}
\email{gmaschler@clarku.edu}
\author[C. T\o nnesen-Friedman]{Christina W.~T\o nnesen-Friedman}
\address{Christina W. T\o nnesen-Friedman\\ Department of Mathematics\\ Union
College\\ Schenectady\\ New York 12308\\ USA } \email{tonnesec@union.edu}

\maketitle
%\tableofcontents
\section{Introduction}
In \cite{guan1}, \cite{guan2}, Guan defined and studied Generalized
Quasi-Einstein (GQE) K\"ahler metrics. On compact manifolds, these are
K\"ahler metrics for which the Ricci potential
is also a Killing potential. This notion includes gradient Ricci solitons
as a special case, and is thus a natural object of study (such solitons
are called Quasi-Einstein metrics in some Physics references).
In \cite{guan2}, GQE metrics are studied in relation to a modified Calabi flow.
Finally, like extremal K\"ahler metrics, GQE metrics generalize the notion of
 constant scalar curvature (CSC) K\"ahler metrics.

Extremal K\"ahler metrics, defined by the requirement that the scalar curvature
is a Killing potential, are the focus of much recent work in K\"ahler geometry.
In \cite{acgt}, a continuity technique was used to show existence of certain
explicit extremal metrics. Our aim in this paper is to apply
the same technique to the question of existence of GQE metrics.

Existence of GQE metrics has been demonstrated in \cite{guan1},
\cite{guan2}, and \cite{ptv} in all K\"ahler
classes of certain manifolds. Here we consider a broader class of spaces,
namely projective bundles over local products of CSC K\"ahler
manifolds that are admissible in the sense defined in \cite{acgt}.
On these spaces we look for a particular type of GQE metric, which we call
admissible. Our main results are as follows. First, we show that any admissible
manifold admits a GQE metric in all admissible K\"ahler classes which are ``small'' in an
appropriate sense. On the other hand, we give an example of a K\"ahler class
on an admissible manifold which is not small, and contains no GQE metric.

%Specifically,
%that the continuity technique used
%to produce extremal K\"ahler metrics on a any admissible manifold (as defined in
%\cite{acgt}) works to produce generalized Quasi-Einstein (GQE)
%K\"ahler metrics as well. Thus, under the assumption that the
%admissible K\"ahler class is ``small'' in an appropriate sense,
%we manage to get a rather general existence result for such metrics.
%We shall see by example that - just like for extremal metrics -
%on certain admissible manifolds there exist K\"ahler classes which are
%not small, having no admissible GQE metrics. For completeness of exposition we shall
%also discuss a part of D. Guan's existence result \cite{guan2} in
%the admissible setting. This existence result - which is proved by
%a root-counting argument - places some scalar curvature sign
%conditions on the base of the admissible
%manifold (albeit not as severe as in the extremal metrics case),
%and in return the smallness condition disappears and
%existence holds in any admissible class on such manifolds.

Our work is laid out as follows.
Section \S\ref{GQE} provides a brief introduction to the Generalized
Quasi-Einstein metrics as defined by D. Guan in \cite{guan1} and
\cite{guan2}.
%These metrics have not yet been explored as thoroughly as
%extremal K\"ahler metrics, but as they specialize to gradient Ricci solitons
%(or quasi-Einstein metrics) when the K\"ahler class is the
%canonical class, they certainly present natural objects of study.
%Furthermore they are explored in connection with a
%modified Calabi flow in, e.g., \cite{guan2}.
Section \S\ref{review} outlines a  brief introduction to the notion of
admissible manifolds, defined in \cite{acgt}, while Section \S\ref{admi}
covers the definition and basic properties of
{\it admissible} Generalized Quasi-Einstein metrics.
Section \S\ref{cont} presents our existence theorem, achieved using a
continuity argument. This is the heart and main purpose of
these notes. Section \S\ref{nonex} provides a non-existence example.
Finally, Section \S\ref{other} contains an appendix discussing another
distinguished metric type of Guan, for which an analog of
the main existence result is obtained.

We would like to thank Vestislav Apostolov for his helpful advice
while preparing this paper.

\section{Background}\label{GQE}

Generalized Quasi-Einstein (GQE) K\"ahler metrics were first introduced by
D. Guan \cite{guan1}. They may be viewed as an alternative
(with respect to extremal K\"ahler metrics) generalization of constant scalar curvature
(CSC) K\"ahler metrics. For instance, any
geometrically ruled surface of genus higher than one has K\"ahler
classes with no extremal metrics (but some K\"ahler
classes on such a manifold do admit extremal metrics \cite{gabor, acgt,t}).
In \cite{ptv} (see also
\cite{guan2} which offers a generalization) it is shown that any
K\"ahler class on this type of manifold admits a GQE metric.

Let $M$ be a complex manifold with almost complex tensor $J$ and a K\"ahler metric
$g$. A function $\phi$ on $M$ is
called a Killing potential if $J\grad\phi$ is a Killing vector field
(i.e. $\nabla J\grad\phi$ is skew-adjoint at every point).
\begin{defn}\cite{guan1,guan2}
\label{def}\label{gqe-def}
Let $g$ be a K\"ahler metric on a compact complex manifold $(M,J)$,
$Scal$  its scalar curvature and $\overline{Scal}$ its average scalar curvature.
We say that $g$ is a GQE metric if there exists a Killing
potential $\phi$ for which
\[
Scal - \overline{Scal} =  \Delta \phi.
\]
Here $\Delta$ denotes the Laplacian with respect to $g$.
\end{defn}
\begin{rem}\label{gqe}
Since $M$ is compact and $Scal - \overline{Scal} = \Delta G \, Scal$,
with $G$ the Green operator, Definition
\ref{def} is equivalent to the requirement that
the Ricci potential  $-G\, Scal$ is also a Killing potential.
In comparison, the definition of
an {\it extremal} K\"ahler metric is equivalent to the statement
that $Scal$ itself is a Killing potential \cite{cal1}.
\end{rem}

\begin{defn}\cite{cal2,fut}
Let $\omega$ be a K\"ahler form on a compact complex manifold $(M,J)$
and let $ h(M)$ denote the Lie algebra of the holomorphic
vector fields on $(M,J)$. Then the Futaki invariant of $[\omega]$ is
the map
$
{\mathcal F}_{[\omega]} : h(M) \rightarrow {\mathbb C}
$
given by
\[
{\mathcal F}_{[\omega]}(\Xi) = - \int_{M} \Xi(G \, Scal) d\mu,
\]
where $\Xi \in h(M)$ and $d\mu$ denotes the volume form
of $\omega$.
\end{defn}

The Futaki invariant is a K\"ahler class invariant. The class of any
 CSC K\"ahler metric has vanishing Futaki invariant. Moreover,
\begin{prop}\cite{guan1}
\label{lucy}
A GQE metric is CSC if and only if the Futaki invariant of the
K\"ahler class vanishes.
\end{prop}

\begin{proof}
We only need to check the ``if'' part of the statement.
Suppose $g$ is a GQE metric as above for some Killing potential $\phi$. Then the
value
$\mathcal{F}_{[\omega]}((\overline{\partial} \phi)^{\sharp})$
of the Futaki invariant on the holomorphic vector field
$(\overline{\partial} \phi)^{\sharp}$ is equal to
\[
-\frac{1}{2}\int_{M}(Scal - \overline{Scal}) \phi \, d\mu
= -\frac{1}{2}\int_{M} \phi \Delta
\phi \, d\mu = - \frac{1}{2}\int_{M} || d\phi||^{2} \, d\mu
\]
(see e.g. \cite{LeSim}).
If this vanishes, then  $\phi$ is constant, and thus
$Scal = \overline{Scal}$.
\end{proof}

\section{Review of admissible manifolds and metrics}\label{review}

Let $S$ be a compact complex manifold admitting a
K\"ahler local product metric, whose components are K\"ahler metrics
denoted $(\pm g_{a}, \pm \omega_{a})$, and
indexed by $a \in {\cA} \subset \mathbb{Z}^{+}$ (Here $\pm g_{a}$
is the K\"ahler metric [In this notation
we allow for the tensors $g_{a}$ to possibly be negatively definite, in which case
the corresponding K\"ahler structure is
$(-g_{a},-\omega_{a})$. If $g_{a}$ is positive definite, then obviously
$(g_{a}, \omega_{a})$ is the corresponding K\"ahler structure.
A parametrization given later justifies this convention.] and $\pm \omega_a$ is the
corresponding K\"ahler form.).
Note that in all our applications, each $\pm g_a$ is assumed to have CSC.
The real dimension of each component is denoted
$2 d_{a}$, while the scalar curvature of $\pm g_{a}$ is given as $\pm 2 d_{a} s_{a}$.
Next, let $E_0$, $E_\infty$ be projectively flat
hermitian holomorphic vector bundles over $S$, of ranks $d_{0}+1$ and $d_{\infty} +1$,
respectively,  such that the condition $c_{1}(E_{\infty})/(d_{\infty}+1) -
c_{1}(E_{0})/(d_{0}+1) = \sum_{a \in\cA} [\omega_{a}/2\pi]$ holds.
Then, following \cite{acgt}, the total space of the projectivization
$M=P(E_{0} \oplus E_{\infty}) \rightarrow S$ is called {\it admissible}.
A particular type of K\"ahler metric on $M$,  also called
{\it admissible}, will now be described.

Let $\hat{\cA} \subset \mathbb{N}\cup\infty$
be the extended index set defined as follows:
\begin{itemize}
\item $\hat{\cA} = \cA$, if $d_{0}=d_{\infty}=0$.
\item $\hat{\cA} =\cA \cup \{ 0 \}$, if $d_{0}>0$ and $d_{\infty}=0$.
\item $\hat{\cA} =\cA \cup \{ \infty \}$, if $d_{0}=0$ and $d_{\infty} >0$.
\item $\hat{\cA} =\cA \cup \{ 0 \} \cup \{ \infty \}$, if $d_{0} > 0$
and $d_{\infty} >0$.
\end{itemize}
In the cases where $\hat{\cA} \neq \cA$, the following notations will prove useful:
$x_{0}=1$, $x_{\infty} = -1$, $s_{0} = d_{0} + 1$ and $s_{\infty} = -(d_{\infty} +1)$.

An admissible K\"ahler metric is constructed as follows. Consider the circle action
on $M$ induced by the standard circle action on $E_0$. It extends to a holomorphic
$\mathbb{C}^*$ action. The open and dense set $M_0$ of stable points with respect to the
latter action has the structure of a principal circle bundle over the stable quotient.
The hermitian norm on the fibers induces via a Legendre transform a function
$z:M_0\ra (-1,1)$ whose extension to $M$ consists of the critical manifolds
$z^{-1}(1)=P(E_{0} \oplus 0)$ and $z^{-1}(-1)=P(0 \oplus E_{\infty})$.
Letting $\theta$ be a connection one form for the Hermitian metric on $M_0$, with curvature
$d\theta = \sum_{a\in\hat{\cA}}\omega_a$, an admissible K\"ahler metric and form are
given up to scale by the respective formulas
\begin{equation}\label{g}
g=\sum_{a\in\hat{\cA}}\frac{1+x_az}{x_a}g_a+\frac {dz^2}
{\Theta (z)}+\Theta (z)\theta^2,\quad
\omega = \sum_{a\in\hat{\cA}}\frac{1+x_az}{x_a}\omega_{a} + dz \wedge
\theta,
\end{equation}
valid on $M_0$. Here $\Theta$ is a smooth function with domain containing
$(-1,1)$ and $x_{a}$, $a \in \cA$ are real numbers of the same sign as
$g_{a}$ and satisfying $0 < |x_a| < 1$. The complex structure yielding this
K\"ahler structure is given by the pullback of the base complex structure
along with the requirement $Jdz = \Theta \theta$. The function $z$ is hamiltonian
with $K= J\grad z$ a Killing vector field, while $\theta$ satisfies
$\theta(K)=1$.

In order that $g$ (be a genuine metric and) extend to all of $M$,
$\Theta$ must satisfy the positivity and boundary
conditions
\begin{align}
\label{positivity}
(i)\ \Theta(z) > 0, \quad -1 < z <1,\quad
%\label{boundary0}
(ii)\ \Theta(\pm 1) = 0,\quad
%\label{boundary1}
(iii)\ \Theta'(\pm 1) = \mp 2.
\end{align}
The last two of these are together necessary and sufficient for
the compactification of $g$.

The K\"ahler class $\Omega_{x} = [\omega]$ of an admissible metric is also called
{\it admissible} and is uniquely determined by the parameters
$x_{a}$, $a \in \cA$, once the data associated with $M$ (i.e.
$d_a$, $s_a$, $g_a$, $z$, $\theta$ etc.) is fixed. The $x_{a}$, $a \in \cA$,
together with the data associated with $M$ will be called {\it admissible data}.
The reader is urged to consult Section 1 of \cite{acgt} for further
background on this set-up.

Define a function $F(z)$ by the formula $\Theta(z)=F(z)/\Mpc(z)$, where
$\Mpc(z) = \prod_{a \in \hat{\cA}} (1 + x_{a} z)^{d_{a}}$.
Since $\Mpc(z)$ is positive for $-1<z<1$, conditions
\Ref{positivity}
imply the following conditions on $F(z)$, which are only necessary
for compactification of the metric $g$:
\begin{align}
\label{positivityF}
(i)\ F(z) > 0, \quad -1 < z <1,\quad
%\label{boundary0F}
(ii)\ F(\pm 1) = 0,\quad
%\label{boundary1F}
(iii)\ F'(\pm 1) = \mp 2p_c(\pm1).
\end{align}

For the purpose of understanding admissible GQE metrics, it is useful to recall the fact below.
\begin{prop}\cite{hamI}
For any admissible metric $g$, if $S(z)$ is a smooth function of z,
then
\begin{equation}
\label{Lapl}
\Delta S = -[F(z) S'(z)]'/p_c(z),\end{equation}
where $\Delta$ is the Laplacian of $g$.
\end{prop}
\begin{proof}
This is a special case of Lemma 5 in \cite{hamI}, but
for convenience we
shall review the proof here.
We denote by $(\_,\_)$ the inner product on two forms induced by
$g$. Recall that
\[
\Delta S = -\left(dd^{c} S(z), \omega\right) = -\left(dJd S(z), \omega\right).
\]
Thus
\[
\begin{array}{ccl}
-\Delta S & = & \left(d (S'(z) J dz\right), \omega)
 =  \left(d (S'(z) \frac{F(z)}{p_{c}(z)} \theta), \omega\right)\\
\\
& = & \left(\left(\frac{[S'(z) F(z)]'}{p_{c}(z)} - \frac{S'(z) F(z)
p_{c}'(z)}{(p_{c}(z))^{2}}\right)dz \wedge \theta, \omega\right)\\
\\
& + & \left(S'(z) \frac{F(z)}{p_{c}(z)} \sum_{a\in\hat{\cA}}\omega_a,
\omega\right)\\
\\
& = & \frac{[S'(z) F(z)]'}{p_{c}(z)} - \frac{S'(z) F(z)}{p_{c}(z)}\left[\frac{
p_{c}'(z)}{(p_{c}(z))}-
\sum_{a\in\hat{\cA}}\frac{d_{a}x_{a}}{(1+x_{a}z)}\right]\\
\\
& = & \frac{[S'(z) F(z)]'}{p_{c}(z)},
\end{array}
\]
where the relation $(\omega_a, \omega) =
(\omega_a, ((1+x_a z)/x_a)\omega_a)
=
( x_a/(1+x_a z))^2 (\omega_a,((1+x_a z)/x_a)\omega_a)_a
=  ( x_a/(1+x_a z)) d_a$, with $(\_,\_)_a$ the inner product
induced by $g_a$, was used.
\end{proof}

The scalar curvature of an admissible metric is given
(cf. \cite{hamI}, or (10) in \cite{acgt})  by
\begin{equation}
Scal =  \sum_{a \in \hat{\cA}}
\frac{2 d_a s_a x_a }{1+x_az} - \frac{F''(z)}{\Mpc(z)},
\label{scalarcurvature}
\end{equation}

Let $C^{\infty}_{*}([-1,1])$ denote the set of functions $f(z)$ of $z$
which are smooth in $[-1,1]$ and normalized so that they integrate to
zero when
viewed as smooth functions on $M$. The latter condition is equivalent
to $\int_{-1}^{1}f(z) p_{c}(z)\, dz =0$, since the volume form of an
admissible metric equals
$p_{c}(z) \big( \bigwedge_{a}
\frac{(\omega_{a}/x_{a})^{d_{a}}}{d_{a}!} \big) \wedge dz \wedge
\theta$.
\begin{cor}
Given an admissible metric $g$, its Laplacian gives a surjective
map from $C^{\infty}_{*}([-1,1])$ to itself (considered as a space
of functions on $M$).
\end{cor}
\begin{proof}
Given $R(z) \in C^{\infty}_{*}([-1,1])$,
an explicit $z$-dependent solution to $\Delta S(z)=R(z)$ can be obtained
directly from \Ref{Lapl} on the open dense set for which for $z\neq\pm 1$.
Either by general Hodge theory or, more concretely, by a L'hospital rule
argument (using (\ref{positivityF}.ii) and
(\ref{positivityF}.iii)), this solution extends to the $\pm 1$ level sets of
$z$.
\end{proof}
\begin{cor}\label{ricpot}
The Ricci potential of an admissible metric is a function of $z$.
\end{cor}
\begin{proof}
This follows from the previous corollary since by
(\ref{scalarcurvature}) the scalar curvature
of an admissible metric is a smooth function of
the moment map $z$.
\end{proof}

\section{GQE metrics on admissible manifolds}\label{admi}

Recall from Remark \ref{gqe} that a K\"ahler metric is GQE if and only if
its Ricci potential is a Killing potential. It follows from Corollary
\ref{ricpot} that an admissible metric $g$ with
moment map $z$ is GQE only if its Ricci potential
is affine in $z$. When this holds, we will call the metric {\it admissible GQE}.
Using Definition \ref{gqe-def} and Remark \ref{gqe},
the admissible GQE condition can be written as
\begin{equation}
\label{gqem1}
Scal - \overline{Scal} =  k \Delta z,
\end{equation}
for some $k \in {\mathbb R}$.

%We have a special case of Proposition \ref{lucy}.

%\begin{prop}
%\label{julia}
%An admissible GQE metric is a CSC metric if and only if
%the Futaki invariant of $(\overline{\partial} z)^{\sharp}$ vanishes.
%\end{prop}

We turn now to an ODE for $F$ which characterizes admissible GQE metrics.
Since for an admissible metric we have from \Ref{scalarcurvature} and \Ref{Lapl} the formulas
\[
Scal =  \sum_{a \in \hat{\cA}}
\frac{2 d_a s_a x_a }{1+x_az} - \frac{F''(z)}{\Mpc(z)},
\quad
\Delta z = \frac{-F'(z)}{\Mpc(z)},
\]
equation \Ref{gqem1} holds if and
only if
\begin{equation}
\label{gqem2}
F''(z) - k F'(z) =
2\biggl(\sum_{a \in \hat{\cA}} \frac{d_a s_a x_a}{1+x_a
z}\biggr)\Mpc(z) - \frac{2\beta_{0}\Mpc(z)}{\alpha_{0}},
\end{equation}

where
\[
\alpha_{0} = \int_{-1}^{1} \Mpc(t)\, dt\quad \mathrm{\ and\ }\quad
 \beta_{0} = \Mpc(1) + \Mpc(-1) +  \int_{-1}^{1}
\biggl(\sum_{a \in \hat{\cA}} \frac{d_a s_a x_a}{1+x_a
t}\biggr)\Mpc(t)\, dt.
\]
Note here that the expression
$\overline{Scal}=2\beta_0/\alpha_0$  (as well as the formula for $Scal$),
appears in the proof Proposition $6$ in \cite{acgt}.

\begin{rem}
Using the extremal polynomial notion (see \cite{acgt}), it is straightforward to verify
that an admissible metric is simultaneously GQE and extremal if and
only if it is CSC. It is tempting to conjecture that this is true in
more generality.
\end{rem}

Just as in the extremal case (see e.g. Section 2.4 in \cite{acgt}), equations
(\ref{positivityF}.ii) and (\ref{positivityF}.iii) together with (\ref{gqem2})
imply (\ref{positivity}.ii) and (\ref{positivity}.iii). So, under
assumption (\ref{gqem2}), (\ref{positivityF}.ii) and (\ref{positivityF}.iii)
are the necessary and sufficient boundary conditions for the compactification
of $g$.

Integrating (\ref{gqem2}) and then solving the resulting first order ODE gives
\begin{equation}
\label{wickfield}
F(z) = e^{k z}\int_{-1}^{z} e^{-k t} P(t)\, dt,
\end{equation}
where $k$ is a constant and
\begin{equation}
\label{Pt}
P(t) = 2 \int_{-1}^{t} \left(\biggl(\sum_{a \in \hat{\cA}} \frac{d_a s_a x_a}{1+x_a
u}\biggr)\Mpc(u) - \frac{\beta_{0}\Mpc(u)}{\alpha_{0}}\right)\, du + 2
\Mpc(-1),
\end{equation}
with the last term determined by the requirement that (\ref{positivityF}.iii)
be satisfied. Also, (\ref{positivityF}.ii) will be satisfied if and
only if there exists
a $k \in \R$ for which
\begin{equation}
\label{kP}
\int_{-1}^{1} e^{-k t} P(t)\, dt = 0.
\end{equation}
In summary, we have
\begin{prop}
\label{lonepine}
Given admissible data on an admissible manifold, let $F$ be the solution of
\Ref{gqem2} of the form \Ref{wickfield}, \Ref{Pt}.
Suppose there exists $k\in \R$ for which \Ref{kP} holds
and (\ref{positivityF}.i) is satisfied by $F$. Then the admissible  metric
naturally constructed from $F$ and the given data is GQE. Conversely,
for any admissible GQE metric (up to scale), the associated function $F$
has the form \Ref{wickfield}, \Ref{Pt}, solves \Ref{gqem2}, satisfies
(\ref{positivityF}.i) and there exists a $k\in R$ for which \Ref{kP} holds.
\end{prop}

We give now two preparatory lemmas on properties of the rational function $P(t)$.
\begin{lemma}
\label{micawber}
For any given admissible data, the function $P(t)$ given by \Ref{Pt}
satisfies: If $d_{0}= 0$, then $P(-1) >0$, otherwise  $P(-1) = 0$.
If $d_{\infty}=0$, then $P(1) < 0$, otherwise $P(1) = 0$.
Furthermore, $P(t)>0$ in some (deleted) right neighborhood of $t=-1$,
and $P(t)<0$ in some (deleted) left neighborhood of $t=1$.
\end{lemma}

\begin{proof}
First observe that by design $P(\pm 1) = \mp 2 \Mpc(\pm 1)$,
which yields the claimed signs of  $P$ at the endpoints.
Also, $\Mpc(t)$ contains the factors $1+x_0t$, $1+x_\infty t$
with multiplicity $d_0$ or, respectively, $d_\infty$.
One of these factors accounts for the vanishing of $P$ at $t=-1$ (or $t=1$)
unless $d_0=0$ (or $d_\infty=0$).
Furthermore, $P'(t)$ contains these factors in each summand, to
multiplicity at least $d_0-1$ (or $d_\infty-1$).
Differentiating P(t), we see that if $d_{0} > 0$, then
$P^{(d_{0})}(-1) > 0$ (and the lower order derivatives vanish),
while if $d_{\infty} > 0$, then $P^{(d_{\infty})}(1)$ has sign $(-1)^{d_\infty +1}$
(and the lower order derivatives vanish).
From these observations the result follows easily by considering the
Taylor expansion of $P(t)$ near $\pm 1$.
\end{proof}

\begin{lemma}
\label{agnes}
If the function $P(t)$ given by \Ref{Pt} has exactly one root in the interval $(-1,1)$, then there
exists a unique $k\in \R$ such that
\begin{equation}\label{k}
\int_{-1}^{1} e^{-k t} P(t)\, dt = 0.
\end{equation}
Moreover, for this $k$, the positivity condition (\ref{positivityF}.i)
is satisfied when $F(z)$ is defined as in (\ref{wickfield}), \Ref{Pt}.
\end{lemma}

\begin{proof}
If $P(t)$ has just one root $t_{0}$ in the interval $(-1,1)$, then,
we may write
\[
P(t) = (t-t_{0})p(t)
\]
where, due to Lemma \ref{micawber}, $p(t)$ is negative for all $t \in (-1,1)$. Consider
now the auxiliary function
\[
G(k) = e^{kt_{0}}\int_{-1}^{1} e^{-k t} P(t)\, dt =
\int_{-1}^{1}p(t)(t-t_{0})e^{-k(t-t_{0})} \, dt.
\]
By direct calculation, $G'(k)$ is positive, while $\lim_{k\rightarrow -\infty}G=-\infty$, and
$\lim_{k\rightarrow \infty}G=+\infty$, as can be checked by taking the limit after
first breaking the integral in the form $\int_{-1}^{t_0}+\int_{t_0}^1$.
This proves the existence and uniqueness
of a $k$ for which $G(k)=0$, or equivalently
$\int_{-1}^{1} e^{-k t} P(t)\, dt = 0$.

Finally, given this $k$, since $e^{-k t} P(t)$ changes sign
exactly once in $(-1,1)$ and is positive near $t=-1$,
condition \Ref{k} clearly guarantees that
$\int_{-1}^{z} e^{-k t} P(t)\, dt$
is a nonnegative function for $z \in (-1,1)$. Therefore
(\ref{positivityF}.i) is satisfied for $F(z)$ as defined in
(\ref{wickfield}), \Ref{Pt}.
\end{proof}

\section{A continuity argument}\label{cont}
Let $M=P(E_{0} \oplus E_{\infty}) \rightarrow S$ be an admissible
manifold, where the base $S$ is a local
K\"ahler product of CSC metrics $(\pm g_{a}, \pm \omega_{a})$.
The aim of this section is
to show that for $|x_{a}|$
sufficiently small for all $a \in \cA$, the corresponding K\"ahler class
admits an admissible GQE metric.
In light of Lemma \ref{agnes}, the strategy will be to
show that in this case $P(t)$ has just one root
in $(-1,1)$.

Observe that
\[
P'(t) = 2\biggl(\sum_{a \in \hat{\cA}} \frac{d_a s_a x_a}{1+x_a
t}\biggr)\Mpc(t) - \frac{2\beta_{0}\Mpc(t)}{\alpha_{0}}
\]
and, as in the proof of Lemma \ref{micawber}, we make the
following observations
\begin{itemize}
\item If $d_{0} > 1$, then $P'(-1) = 0$ and $P'(t)$ is positive in some
(deleted) right
neighborhood of $t=-1$.
\item If $d_0 = 1$, then $P'(-1) > 0$.
\item If $d_{\infty} > 1$, then $P'(1) = 0$ and $P'(t)$ is positive in some
(deleted) left neighborhood of $t=-1$.
\item If $d_{\infty} = 1$, then $P'(1) >0$.
\end{itemize}
We will now look at the behaviour of $P'(t)$ when $x_{a}$ is near $0$ for all
$a \in \cA$. The limit $x_{a} \rightarrow 0$ for all $a\in \cA$ (of any expression)
will be denoted simply by $\lim$. For $P'(t)$, this limit
does not correspond to a K\"ahler class, but is nonetheless a perfectly well
behaved smooth function.
\begin{lemma}
\label{emily}
$\lim P'(t)$, taken as $x_{a} \rightarrow 0$ for all $a \in \cA$, equals
\[
\begin{array}{cl}
  & 2 d_{0}(d_{0}+1)(1+t)^{d_{0}-1}(1-t)^{d_{\infty}}\\
+ & 2 d_{\infty}(d_{\infty}+1)((1+t)^{d_{0}}(1-t)^{d_{\infty}-1} \\
- & (1 + d_{0} + d_{\infty})(2 + d_{0} + d_{\infty})
(1+t)^{d_{0}}(1-t)^{d_{\infty}}.
\end{array}
\]
\end{lemma}

\begin{proof}
The first two terms of the expression simply follows from the fact
%(Theorem 1 in \cite{acgt})
that
$s_{0}x_{0} = d_{0} + 1$ (if $d_{0} \neq 0$) and
$s_{\infty}x_{\infty} = d_{\infty} +1$ (if $d_{\infty} \neq 0$).

The last term follows from the fact that (in the limit considered here)
$\lim (2 \beta_{0}/\alpha_{0})$ equals $(1 + d_{0} + d_{\infty})(2 + d_{0} + d_{\infty})$.
This fact is not at all trivial but follows directly from the calculations
at the end of Appendix B of \cite{acgt}.
\end{proof}

\subsection{Case 1: $d_{0}>0, d_{\infty}>0$}
In this case $\lim P'(t)$ is
\[
g(t)(1+t)^{d_{0}-1}(1-t)^{d_{\infty}-1},
\]
where
\[
g(t) = 2 d_{0}(d_{0}+1)(1-t) + 2 d_{\infty}(d_{\infty}+1)(1+t) -
(1 + d_{0} + d_{\infty})(2 + d_{0} + d_{\infty})(1-t^{2})
\]
is a concave up parabola, which is positive at $t = \pm 1$ and has a
minimum value equal to $-4(1+d_0)(1+d_\infty)/(2+d_0+d_\infty)$,
so negative, in the interval $(-1,1)$.
It is now clear that $\lim P'(t)$ has two distinct {\emph simple} roots
in the interval $(-1,1)$. Thus for
$|x_{a}|$ sufficiently small for all $a\in \cA$, the function $P'(t)$ also has exactly two
zeroes, i.e. $P(t)$ has exactly two [The factored term
$(1+t)^{d_{0}-1}(1-t)^{d_{\infty}-1}$ does not depend on $x_{a}$, so
the corresponding endpoint roots stay put as $x_{a}$ changes.] critical points in
$(-1,1)$.
Putting this together with Lemma \ref{micawber}, we
see that $P(t)$ must change sign exactly once in
$(-1,1)$.

\subsection{Case 2: $d_{0}=0, d_{\infty}>0$}
In this case  $\lim P'(t)$ is
\[
g(t)(1+d_{\infty})(1-t)^{d_{\infty}-1},
\]
where
\[
g(t) =  (2 + d_{\infty}) t + d_{\infty} - 2
\]
is linear and increasing from $g(-1) = -4 < 0$ to $g(1)= 2 d_{\infty}
>0 $. Hence $\lim P'(t)$ has exactly one simple zero
in $(-1,1)$. Thus for
$|x_{a}|$ sufficiently small for all $a\in \cA$, the function $P'(t)$ also has exactly one
zero, i.e. $P(t)$ has exactly one critical point in $(-1,1)$.
Putting this together with Lemma \ref{micawber}, we
see that $P(t)$ must change sign exactly once in
$(-1,1)$.

\subsection{Case 3: $d_{0}>0, d_{\infty}=0$}
In this case $\lim P'(t)$ is
\[
g(t)(1+d_0)(1+t)^{d_0 -1},
\]
where
\[
g(t) =  -(2 + d_0) t + d_0 - 2
\]
is linear and decreasing from $g(-1) = 2 d_{0} > 0$ to $g(1)= -4
>0 $. Hence $\lim P'(t)$ has exactly one simple root
in $(-1,1)$. Thus for
$|x_{a}|$ sufficiently small for all $a\in \cA$, the function $P'(t)$ also has exactly one
zero, i.e. $P(t)$ has exactly one critical point in $(-1,1)$.
Putting this together with Lemma \ref{micawber}, we
see that  $P(t)$ must change sign exactly once in
$(-1,1)$.

\subsection{Case 4: $d_{0}=0=d_{\infty}$}
In this case $\lim P'(t)$ is simply the constant
function $g(t) = -2$. Hence $\lim P'(t)$ has no roots in $(-1,1)$
and is negative. Thus for
$|x_{a}|$ sufficiently small for all $a\in \cA$, the function $P'(t)$ is also strictly
negative, i.e. $P(t)$ is a strictly decreasing function on
 $(-1,1)$. Putting this together with Lemma \ref{micawber}, we
see that $P(t)$ must change sign exactly once in
$(-1,1)$.

Having thus considered all possible cases we may now conclude with
\begin{thm}
\label{copperfield}
Let $M = P(E_{0} \oplus
E_{\infty}) \rightarrow S$ be an admissible manifold arising from a base $S$
with a
local K\"ahler product of CSC metrics. Then the
set of admissible K\"ahler classes admitting an admissible GQE metric
forms a nonempty open subset of the set of all admissible K\"ahler classes.
Any admissible K\"ahler class which is sufficiently small, that is,
for which $|x_{a}|$, $a \in \cA$, are all sufficiently small, belongs to this subset.
\end{thm}

\begin{proof}
The non-emptiness and the inclusion of sufficently small admissible classes
follow from the observations above and Lemma \ref{agnes}.

For the openness we proceed as follows. Recall from Section \ref{review}
that for a given admissible manifold, the admissible K\"ahler classes
are parameterized (up to scale) by $x_{a}$, $a \in \cA$. Suppose
$\cA = \{ 1,\ldots,N \}$, so that the set of admissible
K\"ahler classes (up to scale) is represented by an open subset $W \subset
(-1,1)^{N}$. Rephrasing Proposition \ref{lonepine}, an admissible
K\"ahler class given by $(x_{1},\ldots,x_{N})$ admits an admissible
GQE metric if and only if there exists $k \in {\mathbb R}$ such that
\begin{equation}
\int_{-1}^{1} e^{-kt}P(t) \, dt =0
\label{ui}
\end{equation}
and
\begin{equation}
\int_{-1}^{t} e^{-ku}P(u) \, du > 0, \quad \quad \quad \quad t \in
(-1,1),
\label{lui}
\end{equation}
for $P(t)$ as in \Ref{Pt}.

Suppose that $(x_{1}^{0},\ldots,x_{N}^{0}, k^{0}) \in W \times {\mathbb R}$
satisfies (\ref{ui}) and (\ref{lui}). We need to show that for
$(x_{1},\ldots,x_{N}) \in W$  sufficiently close to
$(x_{1}^{0},\ldots,x_{N}^{0})$, there exists $k \in {\mathbb R}$ such
that $(x_{1},\ldots,x_{N},k)$ also satisfies (\ref{ui}) and (\ref{lui}).
Define $\Phi : W \times {\mathbb R} \rightarrow {\mathbb R}$ by
\[
\Phi(x_{1},\ldots,x_{N},k) = \int_{-1}^{1} e^{-kt}P(t) \, dt,
\]
where $P(t)$ is determined by $(x_{1},\ldots,x_{N})$.
Clearly $\Phi$ is a smooth mapping.
Then
\[
\frac{\partial \Phi}{\partial k} = -\int_{-1}^{1}t e^{-kt}P(t) \, dt
= -\int_{-1}^{1} e^{-kt}P(t) \, dt + \int_{-1}^{1}\big(\int_{-1}^{t}
e^{-ku}P(u) \, du \big) \, dt,
\]
which by (\ref{ui}) and (\ref{lui}) is positive at
$(x_{1}^{0},\ldots,x_{N}^{0}, k^{0})$. A standard implicit function
theorem now gives an open neighborhood $U \subset W$ of
$(x_{1}^{0},\ldots,x_{N}^{0})$ such that for all $(x_{1},\ldots,x_{N}) \in
U$ there exists $k \in {\mathbb R}$ such that
$\Phi(x_{1},\ldots,x_{N},k)=0$, i.e., (\ref{ui}) is satisfied.
Moreover, such $k$ are close to $k^{0}$, when $(x_{1},\ldots,x_{N})$
is close to  $(x_{1}^{0},\ldots,x_{N}^{0})$. By continuity of
$\int_{-1}^{t} e^{-ku}P(u) \, du$ with respect to $x_{1},\ldots,x_{N}$,
and $k$, there is an open neighborhood $V \subset U \subset W$ of
$(x_{1}^{0},\ldots,x_{N}^{0})$ such that for all $(x_{1},\ldots,x_{N}) \in
V$ there exists $k \in {\mathbb R}$ such that (\ref{lui}) as well as
(\ref{ui}) are satisfied. The openness statement now follows and this
concludes the proof of Theorem \ref{copperfield}.
\end{proof}

The theorem above can be compared to D. Guan's existence result
\cite{guan2}, accomplished by a delicate root
counting argument similar in type to the one encountered for extremal K\"ahler
metrics (see Hwang and Singer \cite{hwang-singer} as well as Guan
\cite{guan3}). The argument places some scalar curvature sign
restrictions on the base of the admissible manifold (which are, however,
not as severe as in the corresponding extremal metric case).
Suppose $M$ is a manifold as in Theorem \ref{copperfield} and consider among
the tensors $g_{a}$ two subsets, of
positive definite, and, respectively, negative definite tensors.
Assume that at least one of these subsets has no elements
whose corresponding K\"ahler metric has negative
scalar curvature. Suppose further that the non-zero scalar curvatures
of the K\"ahler metrics corresponding tensors in the complementary subset
all have the same sign [In the most general form of
Cor. 2.13 in \cite{guan2}, the second
condition can be relaxed a bit. Then, however, the existence
appears to depend on the K\"ahler class - in a different sense than our
``smallness'' condition.].
Then Guan's existence result (Cor. 2.13 in \cite{guan2})
implies that
every
admissible K\"ahler class has an admissible GQE metric.

\section{A non-existence example}\label{nonex}
Consider the admissible manifold
\[
P(\cO \oplus \cO(1,-1)) \rightarrow \Sigma_{1} \times \Sigma_{2},
\]
where $\Sigma_{1}$ and $\Sigma_{2}$ are both compact Riemann surfaces
of genus two and $g_{1}$ and $-g_{2}$ are both K\"ahler metrics of
scalar curvature $-4$. Thus $d_{0}=d_{\infty}=0$, $\hat{\cA} = \cA = \{
1,2 \}$, $d_{1}=d_{2}=1$,
$s_{1} = - s_{2} = -2$, and the K\"ahler cone is parametrized by
$0<x_{1}<1$ and $-1 < x_{2} <0$.

Using Proposition 6 in \cite{acgt} one may calculate that the Futaki
invariant of $J\grad z$ equals (up to sign and scale)
\[
\frac{(1+x_{1}-x_{2})(x_{1}+x_{2})}{(3+x_{1} x_{2})^{2}}.
\]

When $x_{2} = -x_{1}$ this vanishes, in fact ${\cF}_{[\omega]}(\Xi)$
vanishes for any $\Xi \in h(M) \cong \C^{\times}$, and using Proposition
\ref{lucy}
we see that any GQE metric in the
corresponding class must be CSC. In turn, any CSC K\"ahler metric must
be admissible \cite{acgt} and thus $k$ in equation
(\ref{gqem1}) should be equal to zero. Calculating $P(t)$ in this
case, we get
\[
P(t) = \frac{2 t (3 - 3x_{1}^{2} - 4 x_{1}^{3} - x_{1}^{2}(1 - 4
x_{1} -x_{1}^{2})t^{2})}{x_{1}^{2} - 3}.
\]
It is easy to see that $\int_{-1}^{1} P(t) \, dt =0$, so
$F(z) = \int_{-1}^{z}P(t) \, dt$ solves (\ref{positivityF}.ii) as well as
(\ref{gqem2})
and (\ref{positivityF}.iii). We calculate that
\[
F(z) = \frac{(1-z^{2})(6 - 7 x_{1}^{2} - 4 x_{1}^{3} + x_{1}^{4}
- x_{1}^{2}(1 - 4
x_{1} -x_{1}^{2})z^{2})}{2(3 - x_{1}^{2})}.
\]
For the interested reader, let us remark that
$F(z)$ is the {\it extremal polynomial} introduced
in \cite{acgt}.

By direct inspection (or by Theorem 2 in \cite{acgt} and Theorem
\ref{copperfield}
in this text), we see
that if $|x_{1}|$ is sufficiently small, (\ref{positivityF}.i) holds
and a CSC metric exists in the corresponding K\"ahler class.
However, for e.g. $x_{1}=0.8$ (and $x_{2}=-0.8$) (\ref{positivityF}.i)
fails, and thus there is exists no GQE metric
in the corresponding K\"ahler class.

Notice, that off but near the line
$x_{2} = -x_{1}$,
(e.g. $x_{1}=
0.9$ and $x_{2}=-0.75$) one may check that there is no extremal
K\"ahler metric in the corresponding class. It can, however,
be shown that in this case $P(t)$ satisfies Lemma \ref{agnes}. Hence this
K\"ahler class {\it admits} an admissible GQE metric.

\begin{rem} It seems to be
``easier'' to
obtain existence of an
admissible GQE metric than that of an (admissible) extremal K\"ahler metric in
a given admissible K\"ahler class.
It is tempting to conjecture that the existence of extremal K\"ahler
metrics in admissible K\"ahler classes
(i.e. positivity of the extremal polynomial) implies the existence of
an admissible GQE metric. Such a result would yield Theorem \ref{copperfield} as
a corollary of Theorem 2 from \cite{acgt}. To determine this one would
have to study more closely the relationship between the extremal
polynomial from \cite{acgt} and $P(z)$.
\end{rem}

\section{Appendix: other metrics}\label{other}
The methods of this paper can be used to give an existence result for
another distinguished metric type, which extrapolates
between extremal and GQE metrics. This type has been considered by Guan
in \cite{guan2}.
Namely, the Killing potential $\phi$ is now
required to satisfy an equation stating that
$\Scal - \overline{Scal}$ is an affine combination of
$\Delta \phi$ and $\phi$. Among admissible metrics
with an associated moment map $z$, we therefore
look for metrics satisfying \begin{equation}
\label{other1}
Scal - \overline{Scal} =  k \Delta z+b(z+l),
\end{equation}
for some $k, b, l \in {\mathbb R}$.
The constant $l$ guarantees that the right hand side of this equation
integrates to zero. It can be computed from admissible data using
its defining equation \Ref{other1}, along with the expressions appearing
in the proof of Proposition $6$ of \cite{acgt}, giving $l=-\alpha_1/\alpha_0$,
with $\alpha_r=\int_{-1}^1p_c(t)t^r \,dt$, $r=1,2$. Using Appendix B of
\cite{acgt}, we have
\begin{lemma}\label{l}
The limit of $l$ as $x_a\rightarrow 0$ for all $a\in \cA$ is
$(d_\infty-d_0)/(2+d_0+d_\infty)$.
\end{lemma}

We now state an existence result for metrics satisfying \Ref{other1}.
\begin{thm}
\label{other thm}
Let $M = P(E_{0} \oplus
E_{\infty}) \rightarrow S$ be an admissible manifold arising from a base $S$
with a local K\"ahler product of CSC metrics.
Then, for any given $b \in {\mathbb R}$, the
set of admissible K\"ahler classes admitting an admissible metric
satisfying \Ref{other1}
forms a nonempty open subset of the set of all admissible K\"ahler classes.
Any admissible K\"ahler class which is sufficiently small, that is,
for which $|x_{a}|$, $a \in \cA$, are all sufficiently small, belongs to this subset.
\end{thm}
\begin{rem} Aside from generalizing Theorem \ref{copperfield}, the above
theorem overlaps with Proposition 9 in \cite{acgt},
which says that for small classes we
may solve \Ref{other1} for $k=0$, obtaining an extremal K\"ahler metric.
Moreover, a solution with $k=0$ can only exist
with a particular - K\"ahler class dependent - value of $b$
(namely $-A$ as defined in Proposition 6 of \cite{acgt},
see also equation (13) there).
Therefore, when $b$ does not equal this
value and is not zero,  Theorem \ref{other thm} guarantees existence
of K\"ahler metrics which are of a new type, i.e. are neither extremal nor
GQE.
\end{rem}
Below we only prove Theorem \ref{other thm} in the case when the
ranks of $E_0$ and $E_\infty$ are at least $2$, i.e. when
$d_{0},d_{\infty}>0$. The general argument is similar.
\begin{proof}
The ODE corresponding to \Ref{gqem2} in this case, is
$$F''(z) - k F'(z) =
2\biggl(\sum_{a \in \hat{\cA}} \frac{d_a s_a x_a}{1+x_a
z}\biggr)\Mpc(z) -
\left(\frac{2\beta_{0}}{\alpha_{0}}+b(z+l)\right)\Mpc(z),$$
and again, assuming this equation holds, (\ref{positivityF}.ii) and
(\ref{positivityF}.iii) are the necessary and sufficient boundary conditions,
which guarantee existence of a metric of type \Ref{other1} on a (compact)
admissible manifold.
Its solution $F$ satisfies, as before,
$F(z) = e^{k z}\int_{-1}^{z} e^{-k t} P(t)\, dt,$
where $P(t)$ (given similarly to (\ref{Pt})) is such that
\[
P'(t) = 2\biggl(\sum_{a \in \hat{\cA}} \frac{d_a s_a x_a}{1+x_a
t}\biggr)\Mpc(t) - \left(\frac{2\beta_{0}}{\alpha_{0}}+b(t+l)\right)\Mpc(t).
\]
For the function $P(t)$, the analog of Lemma \ref{micawber} holds
(since the proof depends largely on $\Mpc(t)$). The analog of Lemma \ref{agnes}
also holds, for fixed $b$ and $l$, with the same proof. Hence what is left is to
analyze $\lim P'(t)$, taken
as
$x_{a} \rightarrow 0$ for all $a \in \cA$. As in Case $1$, we have
$\lim P'(t)=g(t)(1+t)^{d_{0}-1}(1-t)^{d_{\infty}-1}$, yet here $g(t)$ is
the cubic polynomial
\begin{eqnarray*}
g(t) &=& 2 d_{0}(d_{0}+1)(1-t) + 2 d_{\infty}(d_{\infty}+1)(1+t)\\
&-&
(1 + d_{0} + d_{\infty})(2 + d_{0} + d_{\infty})(1-t^{2})-b(t+\lim l)(1-t^2).
\end{eqnarray*}
We have $g(-1)=4d_0(d_0+1)>0$, $g(1)=4d_\infty(d_\infty+1)>0$. Hence
(asymptotics of a cubic show that) one
of the roots of $g(t)$ lies outside $(-1,1)$, and thus at most two lie
in $(-1,1)$. Our proof will be complete once we show that $g(t)$ has exactly
two simple roots in $(-1,1)$, since then the same will hold for $\lim P(t)$,
and we can proceed as in the proof of Theorem \ref{copperfield}.
For this, it is enough to show that $g(t_0)<0$
for some $t_0\in (-1,1)$. Let $t_0=-\lim l=(d_0-d_\infty)/(2+d_0+d_\infty)$.
This number clearly lies in $(-1,1)$, and a direct calculation gives
$g(t_0)=-(4(1+d_\infty)(1+d_0))/(2 + d_{0} + d_{\infty})<0$
%\begin{eqnarray*}
%g(t_0) &=& (1-t_0)\left[2 d_{0}(d_{0}+1) -(1 + d_{0} + d_{\infty})(2 + d_{0} + d_{\infty})
%(1+t_0)\right]\\
%&+&2 d_{\infty}(d_{\infty}+1)(1+t_0)\\
%&=&\frac{2(1+d_\infty)}{2 + d_{0} + d_{\infty}}
%\left[2 d_{0}(d_{0}+1)-(1 + d_{0} + d_{\infty})2(1+d_0)\right]\\
%&+&2 d_{\infty}(d_{\infty}+1)\frac {2(1+d_0)}{2 + d_{0} + d_{\infty}}\\
%&=&\frac{4(1+d_\infty)(1+d_0)}{2 + d_{0} + d_{\infty}}\left[d_0-(1+d_0+d_\infty)+d_\infty\right]\\
%&=&-\frac{4(1+d_\infty)(1+d_0)}{2 + d_{0} + d_{\infty}}<0
%\end{eqnarray*}
as required. This completes the proof of non-emptiness and the inclusion of
sufficiently small admissible classes, using
Lemma \ref{agnes}. Openness follows as in Theorem \ref{copperfield}.
\end{proof}

\begin{rem} It is not hard to check that the K\"ahler class in the
example from Section \ref{nonex}, which carries no GQE nor
extremal K\"ahler metric, does in fact have admissible metrics satisfying
(\ref{other1}).
\end{rem}

%\enddoublespace
\end{document}